\documentclass{conm-p-l}

\newtheorem{theorem}{Theorem}[section]
\newtheorem{lemma}[theorem]{Lemma}
\newtheorem{proposition}[theorem]{Proposition}
\newtheorem{corollary}[theorem]{Corollary}
\newtheorem{conjecture}[theorem]{Conjecture}

\theoremstyle{definition}

\theoremstyle{remark}
\newtheorem{remark}[theorem]{Remark}

\numberwithin{equation}{section}

\newcommand{\BB}{\operatorname{B}}

\newcommand{\HH}{\operatorname{H}}
\newcommand{\CC}{{\mathbb C}}
\newcommand{\Br}{\operatorname{Br}}
\newcommand{\FF}{\operatorname{F}}
\newcommand{\Fp}{{\mathbb F}_p}
\newcommand{\hocolim}{\operatorname{hocolim}}
\newcommand{\map}{\operatorname{Map}}

\begin{document}
\title[Homology of $\Sigma_n$ and polynomial equations]
{A note on the homology of $\Sigma_n$, the Schwartz genus, and
solving polynomial equations}

\author{Gregory Arone}
\address{Department of Mathematics \\
      P.O. Box 400137\\ University of Virginia\\
      Charlottesville VA 22904 - 4137\\ U.S.A.}
\email{zga2m@@virginia.edu}

\thanks{The author was supported in part by
NSF Grant \#0307069}

\date{\today}



\begin{abstract}
We calculate a certain homological obstruction introduced by De
Concini, Procesi and Salvetti in their study of the Schwartz genus
of the fibration $\FF(\CC,n)\to \FF(\CC,n)_{\Sigma_n}$. We show
that their obstruction group vanishes in almost all, but not all,
the hitherto unknown cases. It follows that if $n$ is not a power
of a prime, or twice the power of a prime, then the genus is less
than $n$. The case of $n=2p^k$ where $p$ is an odd prime remains
undecided for some $p$ and $k$.
\end{abstract}
\maketitle
\section{Introduction}
Let $q:E\to B$ be a covering map (or, more generally, a
fibration). By a {\em trivialization} of $q$ one means a
decomposition of $B$ as a finite union of open subsets
$B=\bigcup_{i=1}^t U_i$ such that the restriction of $q$ to
$q^{-1}(U_i)$ is the trivial covering for all $1\le i\le t$. The
{\em Schwartz genus} of $q$ is the minimal $t$ among all such
trivializations of $q$. The Schwartz genus is also known as the
{\em category} of a fibration. In this paper we will refer to the
Schwartz genus as, simply, the genus. Let $g(q)$ be our notation
for the genus of $q$.

This paper is a contribution to the study of the genus of the
quotient map
$$q_n:\FF(\CC,n)\to \FF(\CC,n)_{\Sigma_n}$$ where $\FF(\CC,n)$ is the
configuration space of ordered $n$-tuples of distinct points in
the plane. The interest in the genus of this covering map stems
from its relationship with polynomial equations: The genus of
$q_n$ gives a lower bound on the ``topological complexity'' of any
algorithm for finding the roots of a complex degree-$n$ polynomial
\cite{Smale}. Another way to put it is to say that the Schwartz
genus of $q_n$ gives a lower bound on the number of functions
needed to write the solutions of a polynomial equation of degree
$n$ in terms of the coefficients.

Let us first survey what can be said about the genus of $q_n$ from
general considerations. Suppose $q:E\to B$ is a normal covering
with group $G$. It is well known that the genus of $q$ is less or
equal than $n$ if and only if the classifying map $B\to BG$
factors through the $n-1$-th stage in the Milnor construction for
$BG$, i.e., through the $n-1$-dimensional space $G^{*n}/G$. It
follows that if $B$ is equivalent to an $n-1$-dimensional complex,
then the genus is at most $n$. In our case, it is known that
$\FF(\CC,n)_{\Sigma_n}$ is homotopy equivalent to a CW-complex of
dimension $n-1$ \cite{FN}, and therefore $g(q_n)\le n$. So, the
next interesting question that one may ask about the genus of
$q_n$ is: For which $n$ does $g(q_n)=n$, and for which values of
$n$ does the genus satisfy $g(q_n)<n$?

We build on the recent work of De Concini, Procesi and Salvetti,
who launched an investigation into the cohomological obstruction
for lowering the genus of $q_n$ below the bound given by
dimensional considerations \cite{DCPS}. One of their achievements
was that they succeeded to convert the rather inaccessible
standard cohomological obstruction into a much more tractable
homological obstruction.

We now recall their main result. Let $\Br_n$ be the braid group on
$n$ strings. It is well known that the space
$\FF(\CC,n)_{\Sigma_n}$ can be identified with the classifying
space $B\Br_n$ of the braid group. Furthermore, the map
$\rho_n:B\Br_n\to B\Sigma_n$, induced by the standard group
homomorphism, can serve as a model for the classifying map of
$q_n$. Let $M$ be a module over $\Sigma_n$. We consider $M$ as a
module over $\Br_n$ by pulling back in the obvious way. Standard
obstruction theory implies the following theorem:
\begin{theorem}
$g(q_n)<n$ if and only if the induced homomorphism on cohomology
$$\rho_n^*:\HH^{n-1}(\Sigma_n;M)\to \HH^{n-1}(\Br_n;M)$$
is zero for all $\Sigma_n$-modules $M$.
\end{theorem}
De Concini, Procesi and Salvetti went further than this in that
they showed that it is enough to test that the corresponding
homomorphism on {\em homology} is zero for a certain universal
$\Sigma_n$-module. In more detail, let
$L_n:=\HH^{n-1}(\FF(\CC,n))$ be the top cohomology of the
configuration space. It is well known that $L_n\cong {\mathbb
Z}^{(n-1)!}$. The action of $\Sigma_n$ on $\FF(\CC,n)$ endows
$L_n$ with the structure of a $\Sigma_n$-module. We are interested
in the group homology $\HH_{n-1}(\Sigma_n;L_n)$. It turns out that
this group serves as a home for a universal obstruction for the
Schwartz genus of $q_n$. The following theorem summarizes, and
rephrases slightly, the discussion on page 611 of \cite{DCPS} all
the way up to Theorem 3.2.
\begin{theorem}(De Concini-Procesi-Salvetti) \label{Italians' main theorem}
The Schwartz genus of $q_n$ is less than $n$ if and only if the
induced homomorphism on homology
$$\rho_*:\HH_{n-1}(\Br_n;L_n)\to \HH_{n-1}(\Sigma_n;L_n)$$
is zero.
\end{theorem}
The proof of theorem \ref{Italians' main theorem} utilizes a
specific model for the equivariant topology of $\FF(\CC,n)$
(namely, the authors of \cite{DCPS} use the Salvetti complex). The
theorem does not seem to follow solely from the fact that $L_n$ is
the top cohomology of $\FF(\CC,n)$.

As a corollary, we have the following theorem.
\begin{theorem}\label{crude obstruction} If
$\HH_{n-1}(\Sigma_n;L_n)=\{0\}$ then $g(q_n)<n$.
\end{theorem}

Our goal in this note is to describe the calculation of the
homology groups $\HH_*(\Sigma_n;L_n)$, with special attention to
the group $\HH_{n-1}(\Sigma_n;L_n)$. Our main result (corollary
\ref{mainresult}(1)) is the following:
\begin{theorem}\label{vanishing}
If $n$ is not a power of a prime, or twice the power of a prime,
then $\HH_i(\Sigma_n;L_n)\cong\{0\}$ for all $i$. Therefore, if
$n$ is not a power of a prime or twice the power of a prime, then
$g(q_n)<n$.
\end{theorem}

On the other hand, Vassiliev showed that if $n$ is a prime power
then $g(q_n)=n$ \cite{Vassiliev}. This leaves the case $n=2p^k$
where $p$ is an odd prime. In this case the groups
$\HH_i(\Sigma_{2p^k};L_{2p^k})$ do not vanish for all $i$. Of
course, the question that really interests us is whether the group
$\HH_{2p^k-1}(\Sigma_{2p^k};L_{2p^k})$ vanishes. It turns out that
sometimes it does and some times it does not. We offer the
following partial results to make this point:
\begin{theorem}\label{oddcases}
(1) For all odd primes $p$,
$\HH_{2p-1}(\Sigma_{2p};L_{2p})=\{0\}$. Therefore, if $n=2p$ then $g(q_n)<n$. \\
(2) $\HH_{17}(\Sigma_{18};L_{18})\ne\{0\}$
\end{theorem}
Part (2) of the theorem says that the obstruction group does not
vanish for $n=2\cdot 32$. We are unable to make any conclusions
about the genus of $q_{2p^k}$ for those $p$ and $k$ for which
$\HH_{2p^k-1}(\Sigma_{2p^k};L_{2p^k})\ne\{0\}$.

Theorem \ref{vanishing} is almost implicit in \cite{AM}. However,
it is not made explicit there, and the reader will see that there
are a couple of technicalities to be sorted out. In particular, we
will use a not entirely trivial lemma about Spanier-Whitehead
duality, which is probably of some independent interest (lemma
\ref{Goodwillie}). The lemma has the following little history: I
posted it as a question to Don Davis' discussion list. There were
a few responses, and the best solution (the one that used least
and proved most) was offered by Goodwillie, whose proof we
reproduce in section \ref{mainresults}. It is archived on the
internet \cite{G}.

In any case, the main purpose of this note is to give a concise
description of what is known about $\HH_*(\Sigma_n;L_n)$, and of
what we can conclude about the Schwartz genus at the moment.
Perhaps the most interesting aspect of the paper is the question
that it leaves open: what can one say about the genus of
$q_{2p^k}$ when $\HH_{2p^k-1}(\Sigma_{2p^k};L_{2p^k})$ does not
vanish? In view of Theorem \ref{Italians' main theorem}, this can
be reformulated as a question about differentials in the Serre
spectral sequence for the homology of the fibration sequence
$$\FF(\CC,2p^k)\to \FF(\CC,2p^k)_{\Sigma_{2p^k}}\to\BB\Sigma_{2p^k}$$
taken with coefficients in $L_n$. More precisely, there is a Serre
spectral sequence
$$\HH_i(\Sigma_n;\HH_j(\FF(\CC,n))\otimes L_n) \Rightarrow
\HH_{i+j}(\Br_n;L_n)$$ Taking $n=2p^k$, $i=2p^k-1$, $j=0$, we see
that there is a copy of $\HH_{2p^k-1}(\Sigma_{2p^k};L_{2p^k})$ at
location $(2p^k-1,0)$ of the $E\sp{2}$ term of the spectral
sequence. Theorem \ref{Italians' main theorem} can be rephrased as
saying that $g(q_{2p^k})<2p^k$ if and only if this copy of
$\HH_{2p^k-1}(\Sigma_{2p^k};L_{2p^k})$ gets wiped out by
differentials in the spectral sequence. The author's guess is that
whenever this group is non-zero, it does not get hit by
differentials. In other words, we would like to offer the
following conjecture.
\begin{conjecture}
$g(q_n)=n$ for all $n$ for which $\HH_{n-1}(\Sigma_n;L_n)$ is
non-trivial.
\end{conjecture}

There also is the larger question of determining $g(q_n)$
precisely, rather than just saying whether $g(q_n)<n$. The author
hopes that someone will take up these challenges.

{\em Organization of the paper:} Our starting point for studying
$\HH_*(\Sigma_n;L_n)$ is the identification of $L_n$ with the
homology of a certain familiar space of partitions. We prove this
identification in section \ref{relationtopartitions}, proposition
\ref{model}. We then prove our main results about
$\HH_*(\Sigma_n;L_n)$ in section \ref{mainresults}. The reader
will see that our way to get at $\HH_*(\Sigma_n;L_n)$ is somewhat
roundabout, in that we first consider the homology of $\Sigma_n$
with coefficients in the module $L_n\otimes{\mathbb Z}[-1]$ where
${\mathbb Z}[-1]$ is the sign representation of $\Sigma_n$. Then
we relate $\HH_*(\Sigma_n;L_n)$ with
$\HH_*(\Sigma_n;L_n\otimes{\mathbb Z}[-1])$ using homotopy theory
in a slightly sneaky way. Finally, in section \ref{calculations}
we do some calculations in the non-vanishing case and prove
theorem \ref{oddcases}.

{\em Acknowledgement:} I would like to thank Corrado De Concini,
Claudio Procesi and Mario Salvetti for getting me interested in
the question and for explaining me their work on the subject. I
also thank Tom Goodwillie for the proof of Lemma~\ref{Goodwillie}.

\section{Relation with the poset of
partitions}\label{relationtopartitions} In this section we recall
that $\HH^{n-1}(\FF(\CC,n))$, also known in this paper as $L_n$,
is isomorphic, as a $\Sigma_n$-module, to the homology of the
familiar space of partitions. Let us recall the definition of of
the space of partitions. Let $\underline n=\{1,\ldots,n\}$. Let
$\Lambda$ be the poset (or category) of partitions (equivalence
relations) on $\underline{n}$ ordered by refinements, where we
adopt the convention that $\lambda_1 \le \lambda_2$ if $\lambda_2$
is a refinement of $\lambda_1$. Clearly, $\Lambda$ has an initial
and a final object. Let $\Lambda_i,\Lambda^f$, and $\Lambda_i^f$
be the posets obtained from $\Lambda$ by deleting the initial
object, the final object, and both the initial and final object
respectively. For $n\ge 2$, let $K_n$ be the unreduced suspension
of the geometric realization of $\Lambda_i^f$. $K_n$ is
homeomorphic to the join $S0*|\Lambda_i^f|$. It is well-known that
$K_n$ is homotopy equivalent to a wedge of $(n-1)!$ spheres of
dimension $n-2$ \cite[4.109]{OT}. We define $K_1$ to be the empty
set. We will also need the following generalization of $K_n$: let
$\lambda$ be a partition of $n$ (and assume $\lambda$ is not the
final partition). Consider the category of partitions of $n$ that
are refinements of $\lambda$. Again, this category has an initial
and a final object and we let $K_\lambda$ be the unreduced
suspension of the geometric realization of the category obtained
by removing the initial and the final object. It is easy to see
that $K_\lambda$ is homeomorphic to the join
$K_{\lambda1}*K_{\lambda2}*\cdots*k_{\lambda^i}$ where
$\lambda1,\lambda2,\ldots,\lambda^i$ are the components of
$\lambda$ (so $i=|c(\lambda)|$, where $c(\lambda)$ is the set of
components of $\lambda$). In particular, $K_{\lambda}$ is easily
seen to be equivalent to a wedge sum of spheres of dimension
$n-|c(\lambda)|-1$.

Clearly, the symmetric group $\Sigma_n$ acts on $K_n$, thus making
$\HH_{n-2}(K_n)$ (the only non-trivial reduced homology group of
$K_n$) into a $\Sigma_n$-module. The purpose of this section is to
prove the following proposition:

\begin{proposition}\label{model}
There is an isomorphism of $\Sigma_n$-modules
$$L_n\cong \HH_{n-2}(K_n)$$
\end{proposition}

The proposition is folklore knowledge, but we are not aware of a
precise reference, although it can easily be deduced from various
facts scattered in the literature. Another reason we give our
proof here is that we would like to have a reference to lemma
\ref{filtration} and remark \ref{taylor towers} below.

\begin{proof}[Proof of proposition \ref{model}] Consider $\FF(\CC,n)$ as a subspace $\FF(\CC,n)\subset
\CC^n\subset S^{2n}$. Let $\Delta^nS2$ be the complement of
$\FF(\CC,n)$ in $S^{2n}$. More generally, for a pointed space $X$,
let $\Delta^nX$ be the ``fat diagonal'' in $X^{\wedge n}$. It
follows from duality, together with the fact that $\Sigma_n$ acts
trivially on $\HH_{2n}(S^{2n})$ that the top cohomology of
$\FF(\CC,n)$ is isomorphic, as a $\Sigma_n$-module, to the bottom
homology of $\Delta^nS2$. This bottom homology occurs in dimension
$n$ of $\Delta^nS2$. Thus, there are isomorphisms of
$\Sigma_n$-modules
$$L_n\cong\HH^{n-1}(\FF(\CC,n))\cong \HH_n(\Delta^nS2)$$

The space $\Delta^nS2$, or more generally $\Delta^nX$, can in turn
be thought of as a homotopy colimit over the category $\Lambda^f$.
Indeed, for an object $\lambda$ of $\Lambda$, let $c(\lambda)$ be
the set of components of $\lambda$. It is easy to see that for a
pointed space $X$, diagonal inclusion defines a functor from
$\Lambda$ to spaces given on objects by $\lambda\mapsto X^{\wedge
c(\lambda)}$. It is easy to see that for a well-pointed $X$ there
are equivalences
$$\Delta^nX\cong \mbox{colim}_{\lambda\in\Lambda^f} X^{\wedge
c(\lambda)}\simeq \hocolim_{\lambda\in \Lambda^f}X^{\wedge
c(\lambda)}$$ where we take the colimit and homotopy colimit in
the based category. Taking $X=S2$, we obtain that
$$\Delta^nS2\cong \mbox{colim}_{\lambda\in\Lambda}
S^{2c(\lambda)}\simeq \hocolim_{\lambda\in
\Lambda}S^{2c(\lambda)}$$

We will now introduce a filtration of the functor $\Delta^nX$,
with properties given in the following lemma:
\begin{lemma}\label{filtration}
There exist functors (from spaces to spaces with an action of
$\Sigma_n$) $\Delta^n_iX$ and $\Sigma_n$-equivariant natural
transformation
$$*=\Delta^n_0X\to\Delta^n_1X\to\Delta^n_2X\to\cdots\to\Delta^n_{n-1}X=\Delta^nX$$
where the homotopy cofiber of the map
$\Delta^n_{i-1}X\to\Delta^n_iX$ is $\Sigma_n$-equivalent to
$$\bigvee_{\{\lambda\mid|c(\lambda)|=n-i\}}K_{\lambda}\wedge
X^{\wedge c(\lambda)}$$
\end{lemma}
\begin{remark}\label{taylor towers}
If we apply $\Sigma^\infty$ to the filtration above then the
subquotients become homogeneous functors of $X$. It follows that
the previous lemma is giving a model for the Taylor tower of
$\Sigma^\infty \Delta^nX$. In particular, the $j$-th homogeneous
layer of this functor is
$$\bigvee_{\lambda\mid|c(\lambda)|=j}\Sigma^\infty
K_{\lambda}\wedge X^{\wedge j}$$
\end{remark}

Let $X$ be $k$-connected, $k\ge1$. It follows from the lemma there
is a $\Sigma_n$-equivariant map
$$\Delta^nX\stackrel{\simeq}{\to}\Delta^n_{n-1}
X\to\Delta^n_{n-1}X/\Delta^n_{n-2}X\stackrel{\simeq}{\to}
K_n\wedge X$$ which is $n+2k-1$ connected. In particular, if $X$
is $1$-connected the map is $n+1$-connected. Taking $X=S2$ we see
that there is an $n+1$-connected map $\Delta^nS2\to S2\wedge K_n$.
In particular, it induces an isomorphism on $\pi_n(-)$ which by
Hurewicz theorem is the same as $\HH_n(-)$. Thus, we have an
isomorphism of $\Sigma_n$-modules $L_n\cong\HH_{n-2}(K_n)$. This
completes the proof of proposition \ref{model} \end{proof}
\begin{proof}[Proof of Lemma \ref{filtration}]
We will filter the category $\Lambda$ by the number of components.
For $1\le i\le j\le n$ let $\Lambda_i^j$ be the full subcategory
of $\Lambda$ consisting of partitions $\lambda$ such that $i\le
|c(\lambda)|\le j$. Thus $\Lambda^f=\Lambda_1^{n-1}$, and we have
a system of subcategories
$$\Lambda_{n-1}^{n-1}\hookrightarrow\Lambda_{n-2}^{n-1}\hookrightarrow\cdots\hookrightarrow
\Lambda_1^{n-1}$$ Notice that at each stage the poset is enlarged
by throwing in minimal objects, or to put it differently, there
are morphisms from objects of
$\Lambda_{j}^{n-1}\setminus\Lambda_{j+1}^{n-1}$ to objects of
$\Lambda_{j+1}^{n-1}$, but not the other way around.

We define
$$\Delta^n_iX:=\mbox{colim}_{\lambda\in \Lambda_{n-i}^n}X^{\wedge
c(\lambda)}$$ It is clear that there are $\Sigma_n$-maps
$\Delta^n_{i-1}X \to\Delta^n_iX$ induced by inclusions of
categories. It is not hard to see that in this case colimit is
equivalent to homotopy colimit, so
$\Delta^n_iX\simeq\mbox{hocolim}_{\lambda\in
\Lambda_{n-i}^n}X^{\wedge c(\lambda)}$. It remains to analyze the
homotopy cofiber of the map $\Delta^n_{i-1}X\to\Delta^n_iX$. It is
easy to see that $\Delta^n_{i-1}X$ can be thought of as a homotopy
colimit over the category $\Lambda_{n-i}^n$ rather than
$\Lambda_{n-i+1}^n$ where one extends the functor $\lambda\mapsto
X^{\wedge c(\lambda)}$ to have value $\lambda\mapsto *$ if
$c(\lambda)=n-i$. It follows that the cofiber is equivalent to the
homotopy colimit $\mbox{hocolim}_{\lambda\in\Lambda_{n-i}^n}
G(\lambda)$, where the functor $G:\Lambda_{n-i}^n\to$ Spaces is
defined by
$$G(\lambda)=\left\{\begin{array}{cc} X^{\wedge c(\lambda)} & \mbox{if }
|c(\lambda)|=n-i \\
\ast & \mbox{otherwise} \end{array}\right.$$ It remains to analyze
the homotopy colimit of $G$. Consider again the poset
$\Lambda_{n-i}^n$. Its minimal elements are partitions $\lambda$
which have exactly $n-i$ components. For a partition $\lambda$,
let $\Lambda_{\ge\lambda}$ the the poset of partitions greater
(finer) or equal to $\lambda$. Then the poset $\Lambda_{n-i}^n$
can be thought of as a union
$$\Lambda_{n-i}^n=\bigcup_{\{\lambda\mid|c(\lambda)|=n-i\}}\Lambda_{\ge\lambda}$$
It follows that any homotopy colimit over $\Lambda_{n-i}^n$ can be
written as a homotopy colimit over the category of non-empty
collections of minimal elements of $\Lambda_{n-i}^n$, where to
each such collection $\{\lambda_i\}$ one associates the homotopy
colimit over the intersection $\bigcap_i\Lambda_{\ge\lambda_i}$.
In particular, for our functor $G$, we can write
$$\mbox{hocolim}_{\lambda\in\Lambda_{n-i}^n} G(\lambda)=
\mbox{hocolim}_{\lambda_1,\ldots,\lambda_k}\mbox{hocolim}_{\lambda\in
\bigcap_{i=1}^k \Lambda_{\ge\lambda_i}}G(\lambda)$$ where the
outer homotopy colimit is over finite (non-empty) collections
$\lambda_1,\ldots,\lambda_k$ of distinct minimal elements of
$\Lambda_{n-i}^n$ (i.e., partitions with $n-i$ components). Now
observe that for any such collection $$\bigcap_{i=1}^k
\Lambda_{\ge\lambda_i} =\Lambda_{\ge\bigwedge_{i=1}^k \lambda_i}$$
where by $\bigwedge_{i=1}^k \lambda_i$ we mean the coarsest common
refinement of $\lambda_1,\ldots,\lambda_k$. Note, moreover, that
if $k>1$ then $\bigwedge_{i=1}^k \lambda_i$ has more than $n-i$
components, so $G(\lambda)=*$ for all $\lambda\in\bigwedge_{i=1}^k
\lambda_i$. It follows that whenever $k>1$,
$$\mbox{hocolim}_{\lambda\in \bigwedge_{i=1}^k
\Lambda_{\ge\lambda_i}}G(\lambda)=\mbox{hocolim}_{\lambda\in
\bigwedge_{i=1}^k \Lambda_{\ge\lambda_i}}*=*$$ (because we are
taking pointed homotopy colimit). It follows, finally, that
$$\mbox{hocolim}_{\lambda\in\Lambda_{n-i}^n} G(\lambda)\simeq
\bigvee_{\{\lambda\mid|c(\lambda)|=n-i\}}\hocolim_{\Delta\in
\Lambda_{\ge\lambda}}G(\Delta)$$

It remains to point out that for each minimal $\lambda$,
$\Lambda_{\ge\lambda}$ is a category with an initial object, and
$G\mid\Lambda_{\ge\lambda}$ is a functor which takes value $X$ on
the initial object and the value $*$ on all other objects. It
follows easily that $$\hocolim_{\Delta\in
\Lambda_{\ge\lambda}}G(\Delta)\simeq X\wedge K_\lambda$$ for every
partition $\lambda$ with $n-i$ components. So
$$\mbox{hocolim}_{\lambda\in\Lambda_{n-i}^n} G(\lambda)\simeq
\bigvee_{\{\lambda\mid|c(\lambda)|=n-i\}}X\wedge K_\lambda$$ This
completes the proof of the lemma.
\end{proof}

\section{Proof of the main results}\label{mainresults}
We think of $K_n$ as a topological realization of the module
$L_n$. Since the reduced homology of $K_n$ is concentrated in
dimension $n-2$, it follows from proposition \ref{model} and a
Serre spectral sequence argument that
$$\HH_*(\Sigma_n;L_n)\cong
\HH_{*+n-2}((K_n)_{h\Sigma_n})$$ where by $(K_n)_{h\Sigma_n}$ we
mean the reduced Borel construction:
$(K_n)_{h\Sigma_n}=K_n\wedge_{\Sigma_n} E\Sigma_{n_+}$. In
particular,
$\HH_{n-1}(\Sigma_n;L_n)\cong\HH_{2n-3}((K_n)_{h\Sigma_n})$.

We will want to relate $\HH_*(\Sigma_n;L_n)$ with
$\HH_*(\Sigma_n;L_n\otimes {\mathbb Z}[-1])$ where ${\mathbb
Z}[-1]$ is the sign representation of $\Sigma_n$. Our topological
realization of the module $L_n\otimes{\mathbb Z}[-1]$ is the space
$K_n\wedge S^n$ with the diagonal action of $\Sigma_n$. Obviously,
the homology of this space is concentrated in dimension $2n-2$,
and its non-trivial homology gives the representation
$L_n\otimes{\mathbb Z}[-1]$ since the action of $\Sigma_n$ on
$\HH_n(S^n)$ gives the sign representation. We will also have an
occasion to use the desuspension of this space, which is space
$K_n\wedge S^{n-1}$, where $S^{n-1}$ has an action of $\Sigma_n$
via the standard action on ${\mathbb R}^{n-1}$. This space
realizes the module $L_n\otimes{\mathbb Z}[-1]$ in dimension
$2n-3$.

The following theorem is part of \cite[Theorem 1.1]{AD}:
\begin{theorem}\label{withtwist}
If $n$ is not a power of a prime, then $$(K_n\wedge
S^n)_{h\Sigma_n}\simeq\ast$$ (and therefore also $(K_n\wedge
S^{n-1})_{h\Sigma_n}\simeq\ast$). Moreover, if $n=p^k$ with $k>0$,
then the homology of this space is all $p$-torsion.
\end{theorem}
Next, there is the following theorem, which relates
$(K_n)_{h\Sigma_n}$ and $(K_n\wedge S^n)_{h\Sigma_n}$.

\begin{theorem}\label{SES}
For all $n$, there is a cofibration sequence
$$(K_{\frac{n}{2}}\wedge
S^{\frac{n}{2}})_{h\Sigma_{\frac{n}{2}}}\to (K_n)_{h\Sigma_n}\to
(K_n\wedge S^{n-1})_{h\Sigma_n}$$ (where $K_{\frac{n}{2}}$ is
understood to be a point if $n$ is odd).
\end{theorem}
The following corollary lists some immediate consequences of
theorems \ref{withtwist} and \ref {SES} put together.
\begin{corollary}\label{mainresult}
(1) Unless $n$ is either a power of a prime or twice the power of
a prime, the space $(K_n)_{h\Sigma_n}$ is contractible, and
therefore all the homology groups $\HH_*(\Sigma_n;L_n)$ vanish.

(2) If $n=2^k$, with $k\ge 1$, then there is a long exact sequence
of homology groups
$$\cdots\to\HH_i(\Sigma_{\frac{n}{2}};L_{\frac{n}{2}})\to\HH_i(\Sigma_n;L_n)
\to\HH_{i-n+1}(\Sigma_n;L_n\otimes{\mathbb Z}[-1])\to$$ $$\to
\HH_{i-1}(\Sigma_{\frac{n}{2}};L_{\frac{n}{2}})\to\cdots$$

(3) For every odd prime $p$ there are isomorphisms
$$\HH_*(\Sigma_{p^k};L_{p^k})
\cong\HH_{*-p^k+1}(\Sigma_{p^k};L_{p^k}\otimes{\mathbb Z}[-1])$$
In particular, $\HH_i(\Sigma_{p^k};L_{p^k})\cong\{0\}$ for
$i<p^k-1$.

(4) For every odd prime $p$ there are isomorphisms
$$\HH_*(\Sigma_{p^k};L_{p^k}\otimes {\mathbb Z}[-1])\to\HH_*(\Sigma_{2p^k};L_{2p^k})$$
\end{corollary}

\begin{proof}[Proof of theorem \ref{SES}]
The theorem is almost proved in \cite[Section 4.2]{AM}, but not
quite. What {\em is} proved there is that there is a kind of dual
cofibration sequence
$$[D(K_{\frac{n}{2}}\wedge
S^{\frac{n}{2}})]_{h\Sigma_{\frac{n}{2}}}\leftarrow
[D(K_n)]_{h\Sigma_n}\leftarrow [D(K_n\wedge
S^{n-1})]_{h\Sigma_n}$$ Where $D(-)$ denotes the Spanier-Whitehead
dual of a spectrum. The proof of the dual cofibration sequence
uses calculus of functors and the EHP sequence. Here by the EHP
sequence we mean the sequence of functors $X\to \Omega\Sigma X\to
\Omega\Sigma X^{\wedge 2}$. This sequence is a fibration sequence
in a stable range, and is an actual fibration sequence if $X$ is
an odd-dimensional sphere. Passing to layers in the Goodwillie
tower gives the dual cofibration sequence (details are given in
[op. cit.]). We have a slightly interesting situation here: there
is in fact a sequence of $\Sigma_n$-equivariant maps
$${\Sigma_n}_+\wedge_{\Sigma_{\frac{n}{2}}}D(K_{\frac{n}{2}}\wedge
S^{\frac{n}{2}})\leftarrow D(K_n)\leftarrow D(K_n\wedge S^{n-1})$$
where the composite map is null-homotopic, but which is not
equivalent to a cofibration sequence. However, upon passing to
homotopy orbits, the sequence becomes equivalent to a cofibration
sequence. We would like to conclude that the Spanier-Whitehead
dual sequence has the same property. It is easy to see that the
proof boils down to the following lemma, which may be of some
independent interest:

\begin{lemma}\label{Goodwillie}
Let $G$ be a finite group (actually, $G$ can be any Lie group
whose adjoint representation is orientable) and let $E$ be a
finite spectrum with an action of $G$. Suppose that $E_{hG}\simeq
*$. Then $D(E)_{hG}\simeq *$.
\end{lemma}
\begin{remark}
The above lemma is trivially true if, for instance, $E$ is itself
contractible, or if $E$ is equivalent to a finite spectrum with a
free action of a finite group $G$.
\end{remark}
\begin{remark}
One nice consequence of the lemma is that if a collection of
subgroups of a finite group is {\em ample} in the sense of
\cite[Definition 3.6]{AD}, then, automatically, it is {\em reverse
ample} in the sense of the same definition.
\end{remark}
\begin{proof}[Goodwillie's proof of lemma \ref{Goodwillie}, reproduced from
\cite{G}] The key is to start with the case when $G$ is connected.
If the connected group $G$ acts on the spectrum $E$ and $E$ is
bounded below then the only way $E_{hG}$ can be contractible is if
$E$ is itself contractible (The first nontrivial homotopy group of
$E$ must be the same as that of $E_{hG}$).

Now let $G$ be any compact Lie group and embed it as a subgroup of
a compact connected Lie group $U$. If $G$ acts on $E$ (a spectrum
non-equivariantly equivalent to a finite one), then consider the
induced $U$-spectrum Ind$(E)= U_+ \wedge_G E$. The homotopy orbit
spectrum Ind$(E)_{hU}$ is the same as $E_{hG}$, so the latter is
contractible if and only if Ind$(E)$ is contractible.

The same applies to the $G$-action on the dual spectrum $D(E)$.
And
$$D(\mbox{Ind}(D(E))) = \mbox{Coind}(E) = \map^G(U_+,E)$$
so the question becomes: If the induced $U$-spectrum Ind$(E)$ is
contractible, must the co-induced $U$-spectrum Coind$(E)$ be
contractible, too? The answer is `yes': the two spectra are not
quite the same, but they differ by a little twist that doesn't
matter for our purposes. Namely,
Coind$(E)\simeq\mbox{Ind}(S^{-V}E)$ where $V$ is a representation
of $G$. More precisely, $V$ the difference between the adjoint
representation of $U$ (restricted to $G$) and that of $G$ (this
is, essentially, the Wirthm\"{u}ller isomorphism). And that means
that the homology of Coind$(E)$ and of Ind$(E)$ are related by a
Thom isomorphism as long as the adjoint representation of $G$ is
orientable. So if one is contractible then the other has no
homology and therefore (being bounded below) is contractible.
\end{proof}
This completes the proof of theorem \ref{SES}\end{proof}
\section{Some calculations in the non-vanishing
case}\label{calculations} The purpose of this section is to do
some calculations of
$$\HH_{2p^k-1}(\Sigma_{2p^k};L_{2p^k})=\HH_{4p^k-3}((K_{2p^k})_{h\Sigma_{2p^k}})$$
for $k=1,2$, and to prove theorem \ref{oddcases}.

Let $n=2p^k$, with $p$ and odd prime. By theorems \ref{withtwist}
and \ref{SES},
$$(K_{2p^k})_{h\Sigma_{2p^k}}\simeq (K_{p^k}\wedge
S^{p^k})_{h\Sigma_{p^k}}$$ and so
$$\HH_{4p^k-3}\left((K_{2p^k})_{h\Sigma_{2p^k}}\right)\cong \HH_{4p^k-3}\left((K_{p^k}\wedge
S^{p^k})_{h\Sigma_{p^k}}\right)$$ The homology groups
$\HH_*((K_{p^k}\wedge S^{p^k})_{h\Sigma_{p^k}};\Fp)$ are
calculated in \cite{AM}\footnote{To be precise, the homology
calculated in \cite{AM} is that of the spectrum $(D(K_{p^k})\wedge
S^{p^k})_{h\Sigma_{p^k}}$, but the same method of calcuation
applies to the case that interests us here}. What basically
happens is this: in the standard simplicial model for $K_{p^k}$,
the set of simplices of dimension $k-1$ contains a set isomorphic
to $\Sigma_{p^k}/\wr^k\Sigma_p$, and the homology groups
$\HH_*((K_{p^k}\wedge S^{p^k})_{h\Sigma_{p^k}};\Fp)$ are
isomorphic to the subquotient of
$\HH_*(S^{p^k}_{h\wr^k\Sigma_{p}};\Fp)$ given by the ``completely
inadmissible'' words of length $k$ in the Dyer-Lashoff algebra,
shifted up by $k-1$ degrees.

For instance, in the case $k=1$ $$\HH_*((K_{p}\wedge
S^{p})_{h\Sigma_{p}};\Fp)\cong \HH_*(S^{p}_{h\Sigma_p};\Fp)$$ The
right hand side is, in turn, generated by symbols of the form
$Q^su$ and $\beta Q^su$, $s\ge 1$, where $u$ is of degree $1$ ($u$
is a generator of $\HH_1(S1)$), $Q^s$ is a Dyer-Lashof operation
raising degree by $2s(p-1)$ and $\beta$ is the homology
B\"ockstein, lowering degree by $1$. Thus $\HH_*((K_{p}\wedge
S^{p})_{h\Sigma_{p}};\Fp)$ has a generator $Q^su$ in every
positive dimension that is $1$ modulo $2(p-1)$, has a generator
$\beta Q^s$ in every positive dimension that is $0$ modulo
$2(p-1)$, and is zero otherwise. It follows that the {\em
integral} homology $\HH_*((K_{p}\wedge S^{p})_{h\Sigma_{p}})$ is
non-zero only in dimensions that are $0$ modulo $2(p-1)$. In
particular $\{0\}\cong\HH_{4p-3}((K_{p}\wedge
S^{p})_{h\Sigma_{p}})\cong \HH_{2p-1}(\Sigma_{2p};L_{2p})$. This
proves part (a) of theorem \ref{oddcases}.

Now consider the case $k=2$. The homology groups
$\HH_*((K_{p2}\wedge S^{p2})_{h\Sigma_{p2}};\Fp)$ are generated by
words of the form
$\beta^{\epsilon_1}Q^{s_1}\beta^{\epsilon_2}Q^{s_2}u$ where
$\epsilon_i\in\{0,1\}$, $s_2\ge 1$ and the words are
``inadmissible'' in the Dyer-Lashof algebra, meaning that
$s_1>ps_2-\epsilon_2$. The dimension of such a word is
$1-\epsilon_1+2s_1(p-1)-\epsilon_2+2s_2(p-1)+1=2(s_1+s_2)(p-1)-(\epsilon_1+\epsilon_2)+2$.
In particular, take $p=3$ and consider the word $\beta Q^7Q^1u$.
It gives an element of $\HH_{33}((K_9\wedge
S9)_{h\Sigma_9};{\mathbb F}_3)$. Moreover, this element is in the
image of the B\"ockstein, so it must be a reduction of an integral
element. It follows that $\{0\}\ne\HH_{33}((K_9\wedge
S9)_{h\Sigma_9}\cong\HH_{17}(\Sigma_{18};L_{18})$. This completes
the proof of part (b) of theorem \ref{oddcases}.

\end{document}